\documentclass[twoside,a4paper,reqno,11pt]{amsart} 
\usepackage{amsfonts, amsbsy, amsmath, amssymb, latexsym}
\usepackage{mathrsfs,array}
\usepackage[top=30mm,right=30mm,bottom=30mm,left=30mm]{geometry}
\usepackage{stmaryrd}
\usepackage{bm}

\usepackage{hyperref}

            \def\e{\epsilon} 
             
            \def\ra{\rangle} 
             
            \def\a{\alpha} 
             \def\b{\beta}

            \def\e{\epsilon} 
            \def\O{\Omega} 
             
            \def\l{\lambda}
            \def\L{\langle}
            \def\R{\rangle}

            \def\pf{\noindent {\bf Proof } }

            \def\F{\mathbb{F}}

            \def\d{\delta} 
     \def\D{\Delta}

     \def\g{\gamma} 
      
            \def\hal{\unskip\nobreak\hfil\penalty50\hskip10pt\hbox{}\nobreak 
            \hfill\vrule height 5pt width 6pt depth 1pt\par\vskip 2mm}

           \newtheorem{theorem}{Theorem} 
      
            \newtheorem{thm}{Theorem}[section] 
             
            \newtheorem{lem}[thm]{Lemma} 
             
            \newtheorem{rem}[thm]{Remark}

\begin{document}

\author{Melissa Lee}
 \address{Department of Mathematics,
    Imperial College, London SW7 2BZ, UK}
 \email{m.lee16@imperial.ac.uk}
 
 \author{Martin W. Liebeck}
\email{m.liebeck@imperial.ac.uk}

           \title{Bases of quasisimple  linear groups}
  \begin{abstract} 
Let $V$ be a vector space of dimension $d$ over $\F_q$, a finite field of $q$ elements, and let $G \le GL(V)  \cong GL_d(q)$ be a linear group. A {\it base} of $G$ is a set of
	vectors whose pointwise stabiliser in $G$ is trivial. We prove that if $G$ is a quasisimple group  (i.e. $G$ is perfect and $G/Z(G)$ is simple) acting irreducibly on $V$, then excluding two natural families, $G$ has a base of size at most 6. The two families consist of alternating groups ${\rm Alt}_m$ acting on the natural module of dimension $d = m-1$ or $m-2$, and classical groups with natural module of dimension $d$ over subfields of $\F_q$.
            \end{abstract} 
\date{\today}
\maketitle



\footnotetext{This paper represents part of the PhD work of the first author under the supervision of the second. The first author acknowledges the support of an EPSRC International Doctoral  Scholarship at Imperial College London.}

\section{Introduction}

Let $G$ be a permutation group on a finite set $\O$ of size $n$. A subset of $\O$ is said
to be a {\it base} for $G$ if its pointwise stabilizer in $G$ is trivial. The minimal size of a
base for $G$ is denoted by $b(G)$ (or sometimes $b(G,\O)$ if we wish to emphasize the action). It is easy to see that $|G| \le n^{b(G)}$, so that $b(G) \ge \frac{\log |G|}{\log n}$. A well known conjecture
of Pyber \cite{P} asserts that there is an absolute constant $c$ such that if $G$ is primitive on $\O$, then $b(G) < c\, \frac{\log |G|}{\log n}$.
Following substantial contributions by a number of authors, the conjecture was finally established in \cite{Mar} in the following form: 
there is an absolute constant $C$ such that for every primitive permutation group $G$ of degree $n$,  
\begin{equation}\label{bd}
b(G) < 45\, \frac{\log |G|}{\log n} + C.
\end{equation}
To obtain a more explicit, usable bound, one would like to reduce the multiplicative constant 45 in the above, and also estimate the constant $C$.

Most of the work in \cite{Mar} was concerned with affine groups contained in $AGL(V)$, acting on the set of vectors in a finite vector space $V$ (since the conjecture had already been establish for non-affine groups elsewhere). For these, one needs to bound the base size for a linear group $G \le GL(V)$ that acts irreducibly on $V$.
One source for the undetermined constant $C$ in the bound (\ref{bd}) comes from a 
 key result in this analysis, namely Proposition 2.2 of \cite{LSbase}, in which quasisimple linear groups are handled. This result says that there is a constant $C_0$ such that if $G$ is a quasisimple group acting irreducibly on a finite vector space $V$, then either $b(G) \le C_0$, or $G$ is a classical or alternating group and $V$ is the natural module for $G$; here by the natural module for an alternating group $A_m$ over $\F_{p^e}$ ($p$ prime) we mean the irreducible ``deleted permutation module" of dimension $m-\d(p,m)$, where $\d(p,m)$ is 2 if $p|m$ and is 1 otherwise. This result played a major role in the proof of Pyber's conjecture for primitive linear groups in \cite{LSbase, LSbase2}, which was heavily used in the final completion of the conjecture in \cite{Mar}.

The main result in this paper shows that the constant $C_0$ just mentioned can be taken to be 6. Recall that for a finite group $G$, we denote by $E(G)$ the 
the subgroup generated by all quasisimple subnormal subgroups of $G$. Also write $V_d(q)$ to denote a $d$-dimensional vector space over $\F_q$.

\begin{theorem} \label{main}
     Let $V = V_d(q)$ ($q=p^e$, $p$ prime) and $G \leq GL(V)$, and suppose that $E(G)$ is 
     quasisimple and absolutely irreducible on $V$. Then one of the following 
     holds: 
\begin{itemize}
\item[{\rm (i)}] $E(G) = {\rm Alt}_m$ and $V$ is the natural ${\rm Alt}_m$-module  over $\F_q$, of dimension $d = m-\d(p,m)$; 
\item[{\rm (ii)}] $E(G) = Cl_d(q_0)$, a classical group  with natural module of dimension $d$ over a subfield $\F_{q_0}$ of $\F_q$;
\item[{\rm (iii)}] $b(G)\le 6$. 
\end{itemize}
\end{theorem}

This result has been used in \cite{HLM} to improve the bound (\ref{bd}),  replacing the multiplicative constant 45 by 2, and the constant $C$ by 24. 

With substantially more effort, it should be possible to reduce the constant 6 in part (iii) of the theorem, and work on this by the first author is in progress.

\section{Preliminary lemmas}

If $G$ is a finite classical group with natural module $V$, we define a {\it subspace action} of $G$ to be an action on an orbit of subspaces of $V$, or, in the 
case where $G = Sp_{2m}(q)$ with $q$ even, the action on the cosets of a subgroup $O_{2m}^\pm (q)$. 

\begin{lem}\label{sporex}
Let $G$ be an almost simple group with socle $G_0$, and suppose $G$ acts transitively on a set $\O$.
\begin{itemize}
\item[{\rm (i)}] If $G_0$ is exceptional of Lie type, or sporadic, then $b(G)\le 7$, with equality only if $G = M_{24}$.
\item[{\rm (ii)}] If $G_0$ is classical, and the action of $G$ on $\O$ is primitive and not a subspace action, then $b(G)\le 5$,  with equality if and only if  $G = U_6(2).2$, 
$\O = (G:U_4(3).2^2)$.
\end{itemize}
\end{lem}

\pf Part (i) follows from \cite[Corollary 1]{BLS} and \cite[Corollary 1]{BOW}. Part (ii) is \cite[Theorem 1.1]{Bur}. \hal

For a simple group $G_0$, and $1\ne x \in {\rm Aut}(G_0)$, define $\a(x)$ to be the minimal number of $G_0$-conjugates of $x$ required to generate the group $\langle G_0, x\ra$, and define
\[
\a(G_0) = {\rm max}\,\left(\a(x)\,:\,1\ne x \in {\rm Aut}(G_0) \right).
\]

\begin{lem}\label{gusa} Let $G_0 = Cl_n(q)$, a simple classical group over $\F_q$ with natural module of dimension $n$. Then one of the following holds:
\begin{itemize}
\item[{\rm (i)}] $\a(G_0) \le n$;
\item[{\rm (ii)}] $G_0 = PSp_n(q)$ ($q$ even) and $\a(G_0) \le n+1$;
\item[{\rm (iii)}] $G_0 = L_2(q)$ and $\a(G_0) \le 4$;
\item[{\rm (iv)}] $G_0 = L_3(q)$ and $\a(G_0) \le 4$;
\item[{\rm (v)}] $G_0= L_4^\e(q)$ and $\a(G_0) \le 6$;
\item[{\rm (vi)}] $G_0 = PSp_4(q)$ and $\a(G_0) \le 5$;
\item[{\rm (vii)}] $G_0 = L_2(9),\,U_3(3)$ or $L_4^\e(2)$.
\end{itemize}
\end{lem}

\pf This is \cite[3.1 and 4.1]{GS}. \hal

To state the next result, let $\bar G$ be a simple algebraic group over an algebraically closed field $K$ of characteristic $p$, and let $V=V(\l)$ be an irreducible $K\bar G$-module of $p$-restricted highest weight $\l$. Let $\Phi$ be the root system of $\bar G$, with fundamental roots $\a_1,\ldots ,\a_l$, and let $\l_1,\ldots ,\l_l$ be corresponding fundamental dominant weights. Denote by $\Phi_S$ (resp. $\Phi_L$) the set of short (resp. long) roots in $\Phi$, and if all roots have the same length, just write $\Phi_S=\Phi$, $\Phi_L=\emptyset$. Let $W = W(\Phi)$ be the Weyl group, and for $\a \in \Phi$ let $U_\a = \{u_\a(t):t\in K\}$ be a corresponding root subgroup.

Now let $\mu$ be a dominant weight of $V=V(\l)$, write $\mu = \sum_{j=1}^l c_j\l_j$, and let $\Psi = \langle \a_i : c_i=0\ra$ be a subsystem of $\Phi$. Define
\[
r_\mu = \frac{|W:W(\Psi)|\cdot |\Phi_S\setminus \Psi_S|}{2|\Phi_S|},\;\;
r_\mu' = \frac{|W:W(\Psi)|\cdot |\Phi_L\setminus \Psi_L|}{2|\Phi_L|}
\]
(the latter only if $\Phi_L\ne \emptyset$).  Let
\[
s_\l = \sum_{\mu} r_\mu,\;\; s_\l' = \sum_{\mu} r_\mu'\;(\hbox{ if }\Phi_L\ne \emptyset),
\]
where each sum is over the dominant weights $\mu$ of $V(\l)$.

For $g \in \bar G\setminus Z(\bar G)$ and $\g \in K^*$, let $V_\g(g) = \{v \in V : vg = \g v\}$, and write ${\rm codim}V_\g(g) = \dim V - \dim V_\g(g)$.

\begin{lem}\label{lawth} Let $V = V(\l)$ as above. 
\begin{itemize}
\item[{\rm (i)}] If $g \in \bar G\setminus Z(\bar G)$ is semisimple, and $\g \in K^*$, then ${\rm codim}V_\g(g) \ge s_\l$.
\item[{\rm (ii)}] If $\a \in \Phi_S$, then ${\rm codim}V_1(u_\a(1)) \ge s_\l$.
\item[{\rm (iii)}] If $\Phi_L\ne \emptyset$ and $\b \in \Phi_L$, then ${\rm codim}V_1(u_\b(1)) \ge s_\l'$.
\item[{\rm (iv)}] For any non-identity unipotent element $u \in \bar G$, we have ${\rm codim}V_1(u) \ge {\rm min}(s_\l,s_\l')$.
\end{itemize}
\end{lem}

\pf Parts (i)-(iii) are \cite[Prop. 2.2.1]{law}. For part (iv), note that \cite[Cor. 3.4]{GM} shows that $\dim V_1(u)$ is bounded above by the maximum of $\dim V_1(u_\a(1))$ and $\dim V_1(u_\b(1))$; hence (iv) follows from (ii) and (iii). \hal

For $\bar G$ of type $D_5$ or $D_6$ and $V$ a spin module for $\bar G$, we shall need the following sharper result. Note that the root system $D_n$ has two subsystems of type $A_1^2$ (up to conjugacy in the Weyl group); with the usual labelling of fundamental roots, we denote these by $(A_1^2)^{(1)} = \langle \a_1,\a_3 \rangle$ and 
$(A_1^2)^{(2)} = \langle \a_{n-1},\a_n \rangle$.

\begin{lem}\label{d56}
Let $\bar G = D_n$ with $n\in \{5,6\}$, and let $V=V(\l)$ be a spin module for $\bar G$ with $\l = \l_n$ or $\l_{n-1}$. 
Let $s \in \bar G\setminus Z(\bar G)$ be a semisimple element, and $u \in \bar G$ a unipotent element of order $p$.
\begin{itemize}
\item[{\rm (i)}] Suppose $n=6$. Then ${\rm codim}V_\g(s) \ge 12$ for any $\g \in K^*$; and ${\rm codim}V_1(u) \ge 12$ provided $u$ is not a root element.
\item[{\rm (ii)}] Suppose $n=5$. 
\begin{itemize}
\item[{\rm (a)}] Then ${\rm codim}V_\g(s) \ge 8$ for any $\g\in K^*$, provided $C_{\bar G}(s)'\ne A_4$; and if $C_{\bar G}(s)'= A_4$, then ${\rm codim}V_\g(s) \ge 6$. 
\item[{\rm (b)}] Provided $u$ is not a root element and also does not lie in a subsystem subgroup $(A_1^2)^{(1)}$, we have ${\rm codim}V_1(u) \ge 8$.
\end{itemize}
\end{itemize}
\end{lem}

\pf For semisimple elements $s$, we follow the method of \cite[Section 8]{law} (originally in \cite{ken}). Let $\Psi$ be a subsystem of the root system $\Phi$ of $\bar G$, and define an equivalence relation on the set of weights of $V(\l)$ by saying that two weights are related if their difference is a sum of roots in $\Psi$. Call the equivalence classes $\Psi$-{\it nets}. 

Now define $\Phi_s = \{\a \in \Phi\,|\,\a(s)=1\}$, the root sytem of $C_{\bar G}(s)$. If $\Phi_s \cap \Psi = \emptyset$, then any two weights in a given $\Psi$-net  that differ by a root in $\Psi$ correspond to different eigenspaces for $s$. 

The subsystem $\Phi_s$ is contained in a proper subsystem spanned by a subset of the nodes of the extended Dynkin diagram 
of $\bar G$. Suppose $\Phi_s \ne A_{n-1}$. Then it is straightforward to check that there is a subsystem $\Psi$ that is $W$-conjugate to $(A_1)^{(2)}$ such that $\Phi_s\cap \Psi = \emptyset$. For this $\Psi$ there are $2^{n-2}$ 
$\Psi$-nets of size $2$, and so it follows from the observation in the previous paragraph that  ${\rm codim}V_\g(s) \ge 2^{n-2}$ for any $\g \in K^*$.

Now suppose $\Phi_s = A_{n-1}$. Here there is a subsystem $\Psi$ that is $W$-conjugate to $(A_1^2)^{(1)}$ 
such that $\Phi_s\cap \Psi = \emptyset$. For this $\Psi$ there are $2^{n-5}$ (resp. $2^{n-3}$,  $2^{n-3}$) 
$\Psi$-nets of size 4 (resp. 2,1), and hence ${\rm codim}V_\g(s) \ge 2^{n-4}+2^{n-3}$ for any $\g \in K^*$.
This lower bound is 12 when $n=6$, and 6 when $n=5$. This proves (i) and (ii) for semisimple elements.

Now consider unipotent elements $u \in \bar G$ of order $p$.  Assume first that $p$ is odd.  
Recall that the Jordan form of a unipotent element $u \in D_n$ on the natural module 
determines a partition $\l$ of $2n$ having an even number of parts of each even size; moreover, each such partition corresponds to a single conjugacy class, except when all parts of $\l$ are even, in which case there are two classes, interchanged by a graph automorphism of $D_n$ (see \cite[Chapter 3]{LSbk}). Denote by $u_\l$ (and by $u_\l,u_\l'$ for the exceptional partitions) representatives of the unipotent classes in $\bar G$.
By \cite[\S 4]{Spalt}, if $\mu,\l$ are partitions and $\mu < \l$ in the usual dominance order, then $u_\mu$ lies in the closure of the class $u_\l^{\bar G}$ (or $u_\l'^{\bar G}$).

Suppose $u$ is not a root element, and also is not in a subsystem subgroup $(A_1^2)^{(1)}$ when $n=5$. 
Then it follows from the above that the closure of $u^{\bar G}$ contains $u' = u_\mu$ with $\mu = (3,1^{2n-3})$ or $(2^4,1^{2n-8})$, the latter only if $n=6$. 
Moreover, ${\rm codim} V_1(u) \ge {\rm codim} V_1(u')$ (see the proof of \cite[3.4]{GM}). 
If  $\mu = (3,1^{2n-3})$ , then $u'$  lies in the $B_1$ factor of a subgroup $B_1\times B_{n-2}$ of $\bar G$, and the restriction of $V$ to this subgroup is given by \cite[11.15(ii)]{LSbk}; it follows that $u'$ acts on $V$ with Jordan form $J_2^{2^{n-2}}$, giving the conclusion in this case. And if $\mu = (2^4,1^4)$ with $n=6$, then 
 $u'$  is in $(A_1^2)^{(1)}$, which is contained in a subsystem $A_4$, and the restriction of the spin module $V$ to $A_4$ is given by \cite[11.15(i)]{LSbk}; the lower bound on 
${\rm codim}V_1(u')$ in (i) follows easily from this.

It remains to consider unipotent involutions with $p=2$. The conjugacy classes of these in $\bar G$ are described in \cite[\S 7]{AS} (alternatively in \cite[Chapter 6]{LSbk}). Adopting the notation of \cite{AS}, representatives are $a_l,c_l$ ($l$ even, $2\le l\le n$), and also $a_6'$ in $D_6$ (which is conjugate to $a_6$ under a graph automorphism). These are regular elements of Levi subsystem subgroups $S$, as follows:
\[
\begin{array}{l|ccccccc}
u & a_2&c_2&a_4&c_4&a_6&a_6'&c_6 \\
\hline
S & A_1 & (A_1^2)^{(2)} &  (A_1^2)^{(1)} & A_1 (A_1^2)^{(2)} & (A_1^3)^{(1)} &(A_1^3)^{(2)}& A_1^4 
\end{array}
\]
where $(A_1^3)^{(1)} = \langle \a_1,\a_3,\a_5\rangle$ and $(A_1^3)^{(2)} = \langle \a_1,\a_3,\a_6\rangle$.
The restrictions $V\downarrow S$ can be worked out using \cite[11.15]{LSbk}, from which we calculate $\dim C_V(u)$ for all the representatives:
\[
\begin{array}{r|ccccccc}
u & a_2&c_2&a_4&c_4&a_6&a_6'&c_6 \\
\hline
\dim C_V(u),\,n=5 & 12&8&10&8 & - & - \\
\hline
\dim C_V(u),\,n=6 & 24 & 16& 20 & 16 & 20 & 16 & 16 
\end{array}
\]
The conclusion of the lemma follows. \hal

\section{Bases for some subspace actions}\label{basesub}

Let $G = Cl(V)$ be a simple symplectic, unitary or orthogonal group over $\F_q$, with natural module $V$ of dimension $n$. For $r < n$, denote by ${\mathcal N}_r$ an orbit of $G$ on the set of non-degenerate $r$-subspaces of $V$. The main result of this section gives an upper bound for the base size of the action of $G$ on ${\mathcal N}_r$ when $r$ is very close to $\frac{n}{2}$:

\begin{thm}\label{nrbase}
Let $G_0 = PSp_n(q)\,(n\ge 6)$, $PSU_n(q)\,(n\ge 4)$ or $P\O^\e_n(q)\,(n\ge 7,\,q \hbox{ odd})$, and let $G$ be a group with socle $G_0$ such that $G \le PGL(V)$, where $V$ is the natural module for $G_0$. Define
\[
r = \left\{\begin{array}{l}
\frac{1}{2}\left(n-(n,4)\right), \hbox{ if }G_0 = PSp_n(q), \\
\frac{1}{2}\left(n-(n,2)\right), \hbox{ if }G_0 = PSU_n(q) \hbox{ or }P\O^\e_n(q).
\end{array}
\right.
\]
Then $b(G,{\mathcal N}_r)\le 5$.
\end{thm}


Theorem \ref{nrbase} will follow quickly from the following result. The deduction is given in Section \ref{ded}.

\begin{thm} \label{7/30}
Let $G$ and $r$ be as in Theorem $\ref{nrbase}$, and let $H$ be the stabilizer in $G$ of a non-degenerate $r$-subspace in ${\mathcal N}_r$. 
Let $x \in G$ be an element of prime order. 
Then  one of the following holds:
\begin{itemize}
\item[{\rm (i)}] $\frac{\log |x^G \cap H|}{\log |x^G|}< \frac{1}{2}+\frac{7}{30}$;
\item[{\rm (ii)}]  $G_0 = PSp_8(q)$ and $x$ is a unipotent element with Jordan form $(2,1^6)$.
\end{itemize}
\end{thm}

Our proof is modelled on that of \cite[Thm. 1.1]{tim3}, where a similar conclusion is obtained for the action of $G$ on the set of pairs $\{U, U^\perp \}$ of non-degenerate $n/2$-spaces.

\subsection{Proof of Theorem \ref{7/30}}

We shall give a proof of the theorem just for the case where $G_0$ is a symplectic group $PSp_n(q)$. The proofs for the orthogonal and unitary groups run along entirely similar lines. 

We begin with a lemma on the corresponding algebraic groups. Let $K = \bar \F_q$ and $\bar{G} = PSp_n(K)$, and let $V=V_n(K)$ be the underlying symplectic space. As in Theorem \ref{7/30}, write 
$r=\frac{1}{2}\left(n-(n,4)\right) = \frac{1}{2}n-m$, where $m = \frac{1}{2}(n,4)$. 
Let $\bar H$ be the stabilizer in $\bar G$ of a non-degenerate $r$-subspace, so that 
$\bar{H} = (Sp_{n/2-m}(K) \times Sp_{n/2+m}(K))/\{\pm I\}$.

Write $p = {\rm char}(K)$. When $p=2$, the classes of involutions in $\bar G$ are determined by \cite{AS}: for 
any odd $l\le n/2$, there is one class with Jordan form of type $(2^l,1^{n-2l})$, with representative denoted by $b_l$; and for any nonzero even $l\le n/2$ there are two such classes, with representatives denoted by $a_l,c_l$. These are distinguished by the fact that $(v,va_l)=0$ for all $v\in V$.

\begin{lem}\label{myprop}
With the above notation, if $x$ is an element of prime order in $\bar H$, then  $\dim (x^{\bar{G}} \cap {\bar{H}}) \le N_x$, where $N_x$ is given in  Table $\ref{algtable}$. In the table,  $l_0$ is the multiplicity of the eigenvalue 1 in the action of 
$x$ on $V$, and $a_i$ is the number of Jordan blocks of size $i$ in the Jordan form of $x$.
\end{lem}

\begin{table}[h!]
\centering
\label{algtable}
\begin{tabular}{|c|c|}
\hline 
Type of element $x$ & $N_x$ \\ 
\hline 
semisimple of odd prime order & $\frac{1}{2} \dim x^{\bar{G}} + \frac{1}{4}(n-l_0) +m^2$ \\ 

semisimple involutions & $\left(\frac{1}{2} + \frac{2}{n}\right) \dim x^{\bar{G}}$ \\ 

unipotent of odd prime order & $\frac{1}{2} \dim x^{\bar{G}} + \frac{1}{4}(n- \sum_{i \textrm{ odd}} a_i) +m^2$ \\ 

unipotent involutions of types $b_l$, $c_l$ & $\left( \frac{1}{2}+\frac{2m+1}{n+2}\right) \dim x^{\bar{G}}$  \\ 

unipotent involutions of type $a_l$ & $\left(\frac{1}{2}+\frac{3m}{2n}\right)\dim x^{\bar{G}}$ \\ 
\hline 
\end{tabular} 
\caption{Bounds on $\dim (x^{\bar{G}} \cap {\bar{H}}) $ for elements $x$ of prime order.}
\end{table}

\pf Denote by $V_1$ and $V_2=V_1^\perp$  the $(n/2-m)$- and $(n/2+m)$-dimensional subspaces of $V$ preserved by ${\bar{H}}$.
 First suppose $x \in {\bar{H}}$ is a semisimple element of odd prime order $r$.
Define $\omega$ to be an $r$th root of unity and let $\ell_i$ be the multiplicity of $\omega^i$  $(0 \leq i \leq r-1)$ as an eigenvalue of $x$ in its action on $V$. We further define $y_{ij}$ to be the multiplicity of $\omega^i$ as an eigenvalue of $x$ in its action on $V_j$. Note that $\ell_i = y_{i1} + y_{i2}$. 
Then
\[
\dim x^{\bar{G}} = \frac{n^2+n}{2} - \left(\frac{\ell_0}{2} + \frac{1}{2}\sum_{i=0}^{r-1} \ell_i^2\right),
\]
and furthermore, 
\[
\begin{array}{ll}
\dim (x^{\bar{G}} \cap {\bar{H}}) & = \dim x^{\bar{H}} \leq \frac{n^2+2n}{4}+m^2 - (\frac{1}{2} \ell_0 + \frac{1}{4}\sum_{i=0}^{r-1} \ell_i^2 ) \\
&   = \frac{1}{2}\dim x^{\bar{G}} + \frac{1}{4}(n-\ell_0)+m^2 \\
&  \leq \left( \frac{1}{2}+\frac{1}{n+2}\right) \dim x^{\bar{G}} + m^2.
\end{array}
\]

Now suppose that $x$ is a semisimple involution. Here $C_{\bar G}(x)^0$ is the image modulo ${\pm I}$ of either 
$GL_{n/2}(K)$ or $Sp_l(K) \times Sp_{n-l}(K)$, for some even $l\le n/2$. In the first case, 
$\dim x^{\bar{G}} = n^2/4+n/2$ and so 
\[
\dim (x^{\bar{G}} \cap {\bar{H}}) = \dim x^{\bar{H}} = \frac{1}{2} \dim x^{\bar{G}} +\frac{n}{4} +\frac{m^2}{2} \leq \left(\frac{1}{2}+\frac{1}{n}\right) \dim x^{\bar{G}} + \frac{m^2-1}{2} \leq \left(\frac{1}{2}+\frac{2}{n}\right) \dim x^{\bar{G}}. 
\]
Now consider the second case, where $C_{\bar G}(x)^0 = Sp_l(K) \times Sp_{n-l}(K)$. 
Here $x$ is ${\bar{G}}$-conjugate to $[ -I_l, I_{n-l}]$, and  
$\dim x^{\bar{G}}= nl - l^2 = l(n-l)$. For $j=1,2$, the restriction of $x$ to $V_j$ is $Sp(V_j)$-conjugate to $[ -I_{l_j}, I_{c-l_j}]$ for some even integer $l_j \geq 0$. Noting that $l = l_1+l_2$, we then have
\[
\dim(x^{\bar{G}} \cap {\bar{H}}) = l_1(\frac{n}{2}-m-l_1) + l_2(\frac{n}{2}+m-l_2) \leq \frac{1}{2} \dim x^{\bar{G}} + m(l_2-l_1 ) \leq\left( \frac{1}{2}+\frac{2}{n}\right) \dim x^{\bar{G}}.
\]

Now suppose that $x$ is a unipotent element of odd prime order $r=p$ and that $x$ has Jordan form on $V$ corresponding to the partition $(r^{a_r}, \dots , 1^{a_1}) \vdash n$. We have two further partitions $(r^{b_r}, \dots , 1^{b_1}) \vdash n/2 -m$ and $(r^{c_r}, \dots , 1^{c_1}) \vdash n/2 +m$ associated to $x$ because it preserves $V_1$ and $V_2$. Notice that $a_i = b_i+c_i$.
By \cite[1.10]{LLS}, 
\[
\dim x^{\bar{G}} = \frac{ n^2+n}{2} - \frac{1}{2}\sum_{i=1}^r \Big(\sum_{k=i}^r a_k \Big)^2-\frac{1}{2} \sum_{i \textrm{ odd}} a_i.
\] 
Hence, using  \cite[p.698]{tim3}, we have 
\[
\dim (x^{\bar{G}} \cap {\bar{H}}) \leq \frac{1}{2} \dim x^{\bar{G}} + \frac{1}{4}(n- \sum_{i \textrm{ odd}} a_i) +m^2 \leq \left( \frac{1}{2}+\frac{1}{n+2}\right) \dim x^{\bar{G}} + m^2.
\]

Finally, we consider the case where $x$ is a unipotent involution. First suppose that $x$ is ${\bar{G}}$-conjugate to either $b_l$ or $c_l$ (as described in the preamble to the lemma). Then \cite[1.10]{LLS} implies that $\dim x^{\bar{G}} = l(n-l+1)$. Let $x$ act on $V_i$ with associated partition $(2^{l_i}, 1^{c_i-l_i})$ for $i=1,2$, where $c_1 = n/2-m$ and $c_2=n/2+m$. Then
\[
\dim (x^{\bar{G}} \cap {\bar{H})}= \dim x^{\bar{H}} \leq  \frac{1}{2} \dim x^{\bar{G}} + \frac{l}{2}+m(l_2-l_1) \leq \left( \frac{1}{2}+\frac{2m+1}{n+2}\right) \dim x^{\bar{G}}.
\]
 Lastly, if $x$ is ${\bar{G}}$-conjugate to $a_l$ for some $2\leq l \leq n/2$, then by \cite[1.10]{LLS}, $\dim x^{\bar{G}} = l(n-l)$. By the definition of an $a$-type involution, if $y \in x^{\bar{G}} \cap {\bar{H}}$ fixes a subspace $V_i$, then the restriction of $y$ to $V_i$ is conjugate to $a_{l_i}$ for some even integer $l_i \geq 0$. Therefore
\[
\dim (x^{\bar{G}} \cap {\bar{H}}) = \dim x^{\bar{H}} \leq \frac{1}{2}\dim x^{\bar{G}} +m( l_2-l_1) 
\]
and we determine that $l_2-l_1 < \frac{3l(n-l)}{2n}$, so  
\[
\dim (x^{\bar{G}} \cap {\bar{H}} )= \dim x^{\bar{H}} \leq \left(\frac{1}{2}+\frac{3m}{2n}\right)\dim x^{\bar{G}}.
\]
This completes the proof of the lemma.
\hal

Now we embark on the proof of Theorem \ref{7/30}, considering in turn the various types of elements $x$ of prime order in the symplectic group $G$. We shall frequently use the notaion for such elements given in \cite[\S 3.4]{bg}. Our approach in general is to find a function $\kappa(n)$ such that 
\begin{equation}\label{inkap}
\frac{\log |x^G \cap H|}{\log |x^G|}< \frac{1}{2}+\kappa(n),
\end{equation}
where $\kappa(n)<\frac{7}{30}$ except possibly for some small values of $n$; these small values are then handled separately, usually by direct computation.

\begin{lem}\label{ssodd}
The conclusion of Theorem $\ref{7/30}$ holds when $x$ is a semisimple element of odd order. 
\end{lem}

\pf 
Suppose $x \in H$ is a semisimple element of odd prime order $r$. Let $\mu = (\ell, a_1, \dots, a_k)$ be the tuple associated to $x$ (as defined in \cite[Definition 3.27]{tim2}), and define $i$ to be the smallest natural number such that $r \mid q^i-1$.  According to \cite[3.30]{tim2} this means that 
\[
|C_G(x)| = \left\{\begin{array}{l} |Sp_l(q)|\prod_{j=1}^k |GL_{a_j}(q^i)|,\;i \hbox{ odd} \\
                                                 |Sp_l(q)|\prod_{j=1}^k |GU_{a_j}(q^{i/2})|,\;i \hbox{ even}.
\end{array}
\right.
\]
Let $d$ to be the number of non-zero $a_j$, and further define $e$ to be equal to 1 or 2 when $i$ is even or odd respectively.
By Lemma \ref{myprop} and adapting the argument given in \cite[p.720]{tim3}, we have
\begin{equation}\label{inn1}
|x^G \cap H| <  \left(\frac{n-l}{di}+1\right)^{d/e} 2^{d(e-1)} q^{\tfrac{1}{2} \dim x^{\bar{G}} + \tfrac{1}{4}(n-\ell)+m^2}.
\end{equation}
Furthermore, \cite[3.27]{tim2} implies that
\begin{equation}\label{inn2}
|x^G| \geq \frac{1}{2} \left( \frac{q}{q+1} \right)^{d(2-e)} q^{\dim x^{\bar{G}}}. 
\end{equation}
and \cite[3.33]{tim2} gives the lower bound
\begin{equation}\label{inn3}
\dim x^{\bar{G}} \geq \frac{1}{2}(n^2+n-l^2-l-\frac{1}{ei}(n-l-i(d-e))^2-i(d-e)).
\end{equation}

First suppose $m=1$ (so that $n\equiv 2\hbox{ mod }4$). Then (\ref{inn1})--(\ref{inn3}) imply that the inequality (\ref{inkap}) holds with $\kappa(n) = \frac{3}{n}+\frac{1}{n+1}$. Note that  $\kappa (n)< 7/30$ for $n\geq 18$. For $n=6,10,14$, we must either adjust our value of $\kappa(n)$ or compute $|x^G\cap H|$ and $|x^G|$ explicitly, since 
here $\frac{3}{n} + \frac{1}{n+1}>7/30$. For $n=14$, we find that (\ref{inkap}) holds with $\kappa(n) = 7/30$  for all choices of $(l,i, d)$ except $(l,i, d) = (0,1,2)$. In the latter case, $H =(Sp_8(q)\times Sp_6(q))/\{\pm I\}$ and 
$|C_G(x)| = |GL_{a_1}(q)|\, |GL_{a_2}(q)|$ with $a_1+a_2=7$. Hence 
\[
|x^G\cap H| = \sum_{b_i\le a_i,b_1+b_2=4} |Sp_8(q):GL_{b_1}(q) \times GL_{b_2}(q)| + 
 |Sp_6(q):GL_{a_1-b_1}(q) \times GL_{a_2b_2}(q)|,
\]
and explicit computation gives $\log |x^G \cap H| /\log |x^G| <\frac{1}{2}+\frac{7}{30}$. For $n=10$, (\ref{inkap}) holds with 
$\kappa(n)=7/30$ for all valid choices of $(l,i,d)$ except $(l,i,d)=(0,1,2)$ or $(0,1,4)$, and again explicit calculations as above give $\log |x^G \cap H| /\log |x^G| < \frac{1}{2}+\frac{7}{30}$. Finally, for $n=6$, we find that $\log |x^G \cap H| /\log |x^G| <  \frac{1}{2}+\frac{7}{30}$ for all choices of $x$ with associated parameters $(l,i,d)$.

Now suppose $m=2$.  Then  (\ref{inn1})--(\ref{inn3}) imply that (\ref{inkap}) holds with  $\kappa(n) = \frac{79}{20(n+1)}$ (when $e=1$), and with $\kappa(n) =\frac{22}{5(n+2)}$ (when $e=2$). We have $\kappa(n) < \frac{7}{30}$ for $n\ge 20$. For $n<20$, explicit calculations of $|x^G\cap H|$ as above yield the conclusion.
\hal

\begin{lem}\label{ssinvol}
The conclusion of Theorem $\ref{7/30}$ holds when $x$ is a semisimple involution. 
\end{lem}

\pf  
Suppose that $x \in H$ is a semisimple involution. Denote by $s$ the codimension of the largest eigenspace of $x$ on 
$V = V_n(K)$. According to \cite[3.37]{tim2}, $|C_G(x)|$ is equal to $|Sp_s(q)|\,|Sp_{n-s}(q)|$, 
$|Sp_{n/2}(q)|^2.2$, $|Sp_{n/2}(q^2)|.2$ or $|GL_{n/2}^\e(q)|.2$, with $s<\frac{n}{2}$ in the first case, and $s=\frac{n}{2}$ in the latter three cases. 
Suppose $x$ is as in one of the first two cases. Adapting the analogous argument given in \cite[p.720]{tim3}, we deduce that
\[
|x^G \cap H | < 4\left(\frac{q^2+1}{q^2-1} \right) q^{\frac{s(n-s)}{2} - m(1-m)},\;\; |x^G| > \frac{1}{2}q^{s(n-s)}
\]
(the constant $\frac{1}{2}$ in the second inequality should be replaced by $\frac{1}{4}$ when $s=\frac{n}{2}$).
These bounds imply that (\ref{inkap}) holds with 
\[
\kappa(n) = \left\{\begin{array}{l} \frac{2}{n}, \hbox{ if }s<\frac{n}{2}, m=1 \\
\frac{3}{n+1}, \hbox{ if }s<\frac{n}{2}, m=2 \\
\frac{3}{2n}, \hbox{ if }s=\frac{n}{2}, n\ge 12. 
\end{array}
\right.
\]
For $n\ge 12$ we have $\kappa(n)< \frac{7}{30}$, giving the conclusion. And for smaller values of $n$, we obtain the conclusion by explicit calculation of the values of $|x^G\cap H|$ and $|x^G|$.

Next suppose $|C_G(x)| = |Sp_{n/2}(q^2)|.2$. Then $|x^G| > \frac{1}{4}q^{n^2/4}$ by \cite[3.37]{tim2}. If $\frac{n}{4}$ is even then $x^G\cap H= \emptyset$, so assume $\frac{n}{4}$ is odd. An argument analogous to that at the top of p.722 of \cite{tim3} for this case gives $|x^G\cap H| < \frac{1}{4}q^{(n^2/8)+2}$. These bounds imply that (\ref{inkap}) holds with $\kappa(n) = \frac{2}{n}$, and this is less than $\frac{7}{30}$ for all $n\ge 12$.

Finally, suppose that $|C_G(x)| = |GL_{n/2}^\e(q)|.2$. Again \cite[3.37]{tim2} and arguments of \cite[p.722]{tim3} give 
\[
|x^G| > \frac{1}{4}\left(\frac{q}{q+1}\right) q^{\frac{1}{4}n(n+2)},\;\;
|x^G\cap H| < \frac{1}{4} q^{\frac{n^2}{8}+\frac{n}{2}+\frac{m^2}{2}}.
\]
Hence (\ref{inkap}) holds with $\kappa(n) = \frac{5}{2n}$, which is less than $\frac{7}{30}$ for $n>10$, and for $n\le 10$ we obtain the conclusion as usual by explicit calculation of $|x^G\cap H|$ and $|x^G|$.
\hal

\begin{lem}\label{unipodd}
The conclusion of Theorem $\ref{7/30}$ holds when $x$ is a unipotent element of odd order. 
\end{lem}

\pf 
Let $x \in H$ be a unipotent element of order $p$, and suppose $p$ is odd. Let the Jordan form of $x$ on $V$ correspond to the partition $\lambda \vdash n$. By Lemma \ref{myprop},
\begin{equation}
\label{unidim}
\dim{x^{\bar{H}}} \leq \frac{1}{2} \dim x^{\bar{G}} + \frac{1}{4}(n-e)+m^2,
\end{equation}
where $e$ is the number of odd parts in $\lambda$.

\noindent \textbf{Case $\lambda = (k^{n/k})$}

Since $k$ must divide both $n/2-m$ and $n/2+m$, we have $k=2$ or 4 (the latter only if $m=2$). Arguing as at the bottom of p.722 of \cite{tim3}, we have $\dim x^{\bar{G}} \geq \tfrac{1}{4}n(n+2)$, and also
\[
|x^G|  > \frac{q}{q+1} q^{\dim x^{\bar{G}}},\;\;
|x^G \cap H| = |x^H| < 4 q^{\dim x^{\bar{H}}} \leq 4 q^{\frac{1}{2} \dim x^{\bar{G}}+ \frac{1}{4}(n-e)+m^2}.
\]
These bounds imply that (\ref{inkap}) holds with $\kappa(n)=\frac{3}{n+1}$, which is less than $\frac{7}{30}$ for $n\ge 14$. As usual, for smaller values of $n$ we obtain the result by explicit computation of $|x^G\cap H|$ and $|x^G|$.

\vspace{2mm}
\noindent \textbf{Case $\lambda = (2^j, 1^{n-2j})$, $n-2j>0$}

First suppose that $j=1$. Then $|x^G| > \frac{1}{4} q^n$ and $|x^G \cap H |  <   q^{n/2 +m} + q^{n/2-m}$. This implies that $\frac{\log |x^G \cap H|}{\log |x^G|}< \frac{1}{2}+\frac{7}{30}$ for all values of $n\ge 6$ except $n=8$. The case $n=8$ is the exception in part (ii) of Theorem \ref{7/30}.

Next suppose that $j=2$. Here $|x^G| > \frac{1}{4(q+1)} q^{2n-1}$. Since the two Jordan blocks of size 2 can lie in the  two different subspaces $V_1$ and $V_2$, or in the same one, we have 
\[
|x^G \cap H| < q^{(n-2m)/2 + (n+2m)/2} + 2q^{n-4+m(m-1)} + 2q^{n+m(m-1)}
\]
Hence (\ref{inkap}) holds with $\kappa(n)=\frac{3}{n+1}$, which is less than $\frac{7}{30}$ for $n\ge 12$. For smaller values of $n$ we obtain the conclusion by explicit computations of $|x^G\cap H|$ and $|x^G|$.

Finally, assume $j \geq 3$ (and so $n\geq 8$ since $n-2j>0$). The number of ways to distribute the $j$ Jordan blocks of size 2 amongst the subspaces $V_1,V_2$ is at most $j+1$. Then, adapting the analogous bound in \cite[p.723]{tim3} and making use of Lemma \ref{myprop}, we have 
\[
|x^G \cap H|  < 4(j+1)  q^{ \dim x^{\bar{G}}/2+ j/2+m^2}
\]
and as in \cite[p.723]{tim3}, we have  
$|x^G|> \frac{1}{4} q^{\dim x^{\bar{G}}} = \frac{1}{4} q^{j(n-j+1)}$.
 This yields (\ref{inkap}) with  $\kappa(n) = \frac{4}{n+2}$, which is less than $\frac{7}{30}$ for $n\ge 16$.
As usual, smaller values of $n$ are handled by direct computation. 

\vspace{2mm}
\noindent \textbf{Case $\lambda = (k^{a_k}, \dots ,  2^{a_2}, 1^l)$, $k \le n/2+m$ }

In the computations below, we adapt the arguments on p.723 of \cite{tim3}. Let $d$ be the number of non-zero $a_i$. 
Then 
\[
|x^G| > \frac{1}{2^{d+1}} \left( \frac{q}{q+1}\right)^d q^{\dim x^{\bar{G}}}.
\]

If $d=1$ then $\lambda = (k^{(n-l)/k}, 1^l)$, and we can take $k>2$ by the previous case. 
By \cite[1.10]{LLS}, we have 
\[
\dim x^{\bar{G} }= \frac{n^2}{2}+ \frac{n}{2} - \frac{l(n-l)}{k}-\frac{l^2}{2}-\frac{1}{2k}(n-l)^2-\frac{l}{2}-\frac{\alpha}{2k}(n-l),
\]
where $\alpha$ is zero if $k$ is even and one if $k$ is odd.
Arguing as in \cite[p.723]{tim3} we also have
\[
|x^G \cap H| < \left(\frac{n-l}{k} +1\right) 2^2 q^{\dim x^{\bar{G}}/2+ (n-l)(1-\alpha/k)/4 +m^2}.
\]
These bounds imply (\ref{inkap}) with $\kappa(n) = \frac{3}{n-3}$, which is less than $\frac{7}{30}$  for $n\ge 16$, and smaller values of $n$ are handed by explicit computation.

Now suppose that $d \geq 2$.
By \cite[p. 723]{tim3}, 
\[
\dim x^{\bar{G}}\geq \frac{1}{4}n^2 + \frac{1}{4}(d^2-d+2) - \frac{1}{16}d^4-\frac{1}{24} d^3 + \frac{3}{16}d^2-\frac{1}{3}d-\frac{1}{4}l^2-\frac{1}{2}, 
\]
and adapting the analogous bound given in \cite[p.723]{tim3} and referring to Lemma \ref{myprop}, we have
\[
|x^G\cap H| < 4^d \left(\frac{n/2-d^2/4+d/4-l/2-1}{d} +1 \right)^d q^{\tfrac{1}{2}\dim x^{\bar{G}} + (n-l)/4+m^2}
\]
These bounds give (\ref{inkap}) with  $\kappa(n)=\frac{4}{n}$, 
which is less than $\frac{7}{30}$  for $n\ge 18$, and smaller values of $n$ are handed by explicit computation.
\hal

\begin{lem}\label{unipodd}
The conclusion of Theorem $\ref{7/30}$ holds when $x$ is a unipotent involution. 
\end{lem}

\pf Let $p=2$, and recall the description of the involution class representatives $a_l,b_l,c_l$ of $G$ in the preamble to Lemma \ref{myprop}. 

First assume that $x$ is conjugate to $a_l$ for some even integer $l$ with $2\leq l \leq n/2$. If $l=2$, then by \cite[1.10]{LLS} and \cite[Proposition 3.9]{tim2} we have
\begin{equation}
\label{l=2}
|x^G \cap H| < 2 q^{2(n/2-m-2)} + 2q^{2(n/2+m-2)}.
\end{equation}
If $l\geq 4$ then we may adapt the analogous equation in \cite[p.723]{tim3} and obtain
\[
|x^G \cap H| <  (\frac{l}{2}+1) 2^2 q^{(\frac{1}{2} + \frac{3m}{2n}) l(n-l)}.
\]
Furthermore, for all $l$, by \cite[p.723]{tim3}
\[
|x^G| > \frac{1}{2} q^{l(n-l)}.
\]
These bounds imply that $\frac{\log |x^G \cap H|}{\log |x^G|}< \frac{1}{2}+\frac{7}{30}$, provided $n\ge 14$ when $l=2$, and $n\ge 24$ when $l\ge 4$. Smaller values of $n$ can be dealt with by explicit computation of $|x^G\cap H|$ and $|x^G|$.

Now suppose that $x$ is conjugate to either a $b_l$- or $c_l$-type involution. If $l=1$ then by \cite[1.10]{LLS} and \cite[Proposition 3.9]{tim2} 
\begin{equation}
|x^G \cap H| < q^{n/2-m}+ q^{n/2+m},
\label{l=1}
\end{equation}
and if $l=2$, then 
\begin{equation}
\label{l=22}
|x^G\cap H| < q^n + q^{2(n/2-m-1)}+ q^{2(n/2+m-1)}.
\end{equation}
If $l\geq 3$, then by adapting the analogous argument in \cite[p.724]{tim3}, we deduce
\[
 |x^G \cap H|< 4\left(\frac{q^2+1}{q^2-1} \right) (q^{\tfrac{1}{2}\dim x^{\bar{G}}+2m-1}+q^{\tfrac{1}{2}\dim x^{\bar{G}}+m-1})   +4 \left( \frac{q^2+1}{q^2-1}\right) q^{\tfrac{1}{2} \dim x^{\bar{G}} + l/2+ m}
\]
where $\mathrm{dim} x^{\bar{G}} = l(n-l+1)$.
Lastly, \cite[p. 724]{tim3} gives
\[
|x^G| > \frac{1}{2} q^{l(n-l+1)} .
\]
As usual, these bounds imply that $\frac{\log |x^G \cap H|}{\log |x^G|}< \frac{1}{2}+\frac{7}{30}$ for $n\ge 14$, and explicit computations give the same conclusion for smaller values of $n$.
\hal

This completes the proof of Theorem \ref{7/30}.

\subsection{Deduction of Theorem \ref{nrbase}}\label{ded}

The deduction of Theorem \ref{nrbase} from Theorem \ref{7/30} proceeds along the lines of the proof of \cite[1.1]{Bur}.

First we shall require a small extension of \cite[Prop. 2.2]{Bur}.
For a finite group $G$, define 
\[
\eta_G(t) = \sum_{C\in \mathcal{C}} |C|^{-t}
\]
where $\mathcal{C}$ is the set of conjugacy classes of elements of prime order in $G$. 

\begin{lem}\label{zeta} Let $G$ be a finite classical group as in Theorem $\ref{nrbase}$, with $n\ge 6$. 
\begin{itemize}
\item[{\rm (i)}] Then $\eta_G(\frac{1}{3}) < 1$.
\item[{\rm (ii)}] Let $G = PGSp_8(q)$. Then $\eta_G(\frac{1}{3}) < 0.396$.
\end{itemize}
\end{lem}

\pf 
(i) This is \cite[Prop. 2.2]{Bur}.

(ii) We compute the sizes of the conjugacy classes with each centraliser type using \cite[Table B.7]{bg}, and bound the number of classes with each centraliser type using the same arguments as those given in the proof of \cite[Lemma 3.2]{Bur}. The result follows from these computations. \hal

We also need to cover separately the two cases of Theorem \ref{nrbase} for dimensions less than 6.

\begin{lem}\label{u4u5} Theorem $\ref{nrbase}$ holds for $G_0 = PSU_4(q)$ or $PSU_5(q)$.
\end{lem}

\pf Consider the first case, Here $G = PGU_4(q)$ acting on ${\mathcal N}_1$, the set of non-degenerate 1-spaces. Let 
$v_1,\ldots ,v_4$ be an orthonormal basis of the natural module for $G$. If $q$ is odd, then $\L v_1\R$, $\L v_2\R$, 
$\L v_3\R$, $\L v_1+v_2+v_3+v_4\R$ is a base for the action of $G$; and if $q$ is even, then 
 $\L v_1\R$, $\L v_2\R$, $\L v_3\R$, $\L v_1+v_2+v_3 \R$, $\L v_2+v_3+v_4\R$ is a base.

Now let $G = PGU_5(q)$  acting on ${\mathcal N}_2$. Let 
$v_1,\ldots ,v_5$ be an orthonormal basis. Any element of $G$ that fixes the three non-degenerate 2-spaces $\L v_1,v_2\R$, 
$\L v_2,v_3\R$ and $\L v_3,v_4\R$ also fixes $\L v_1,v_5\R$ and $\L v_4,v_5\R$ (as these are $\L v_2,v_3, v_4\R^\perp$ and 
$\L v_1,v_2, v_3\R^\perp$), hence fixes all the 1-spaces $\L v_1\R,\ldots ,\L v_5\R$. Hence adding two further non-degenerate 2-spaces intersecting in $\L v_1+ \cdots +v_5\R$ to the first three gives a base of size 5.
\hal

\vspace{2mm} {\it Proof of Theorem \ref{nrbase} }
Let $G,r$ be as in the statement of Theorem \ref{nrbase}, and let $H$ be the stabilizer of a non-degenerate $r$-subspace in ${\mathcal N}_r$. In view of Lemma \ref{u4u5}, we may assume that the dimension $n\ge 6$. 
 
For a positive integer $c$, let $Q(G,c)$ be the probability that a randomly chosen $c$-tuple of elements of ${\mathcal N}_r$ does not form a base for $G$. Then
\begin{equation}\label{qgc}
Q(G,c) \leq \sum_{x \in X} |x^G| \left( \frac{{\rm fix}_{{\mathcal N}_r}(x)}{|{\mathcal N}_r|}\right)^c =   \sum_{x \in X} |x^G| \left( \frac{|x^G\cap H|}{|x^G|}\right)^c,
\end{equation}
where $X$ is a set of conjugacy class representatives of the elements of $G$ of prime order. Clearly $G$ has a base of size $c$ if and only if $Q(G,c)<1$. 

Assume for the moment that $G_0 \ne PSp_8(q)$. Then by Theorem \ref{7/30} we have 
\[
\frac{|x^G\cap H|}{|x^G|} < |x^G|^{-\frac{1}{2}+\frac{7}{30}}
\]
for all elements $x \in G$ of prime order. Hence it follows from (\ref{qgc}) that 
\[
Q(G,5) < \sum_{x \in X} |x^G|^{1 +5\left(-\frac{1}{2}+\frac{7}{30}\right)} = \eta_G(1/3).
\]
Therefore by Lemma \ref{zeta}(i), $G$ has a base of size 5, as required.

 It remains to consider the case where $G_0 = PSp_8(q)$. Here Theorem \ref{7/30}(ii) gives $\frac{|x^G\cap H|}{|x^G|}< |x^G|^{-\frac{1}{2}+\frac{7}{30}}$ for all elements $x\in G$ of prime order, except when $x$ is a unipotent element with Jordan form 
$(2, 1^6)$. In the latter case $|x^G| = q^8-1$ and $|x^G\cap H| = q^6+q^2-2$.
Hence 
\[
Q(G,5) <  \eta_G(1/3) + (q^8-1)\left(\frac{q^6+q^2-2}{q^8-1}\right)^{5},
\]
and this is less than 1 for all $q$, by Lemma \ref{zeta}(ii).

\vspace{2mm}
This completes the proof of Theorem \ref{nrbase}.


\section{Proof of Theorem \ref{main}}

Assume the hypotheses of Theorem \ref{main}. Thus $G \le GL(V) = GL_d(q)$, and $E(G)$ is quasisimple and absolutely irreducible on $V$. Then the group $Z:=Z(G)$ consists of scalars, and $G/Z$ is almost simple. Let $G_0$ be the socle of $G/Z$. 

\begin{lem}\label{sporadexcep}
If $G_0$ is exceptional of Lie type or sporadic, then $b(G) \le 6$.
\end{lem}

\pf Pick $v \in V\setminus 0$, and consider the action of $G$ on the orbit $\D = v^G$. By Lemma \ref{sporex}(i), if $G_0 \ne M_{24}$ then there exist $Z$-orbits $\d_1,\ldots ,\d_6$ such that $G_{\d_1\cdots \d_6} \le Z$. Hence $b(G) \le 6$. The case where $G_0=M_{24}$ is taken care of in Remark \ref{m24} below. \hal

\begin{lem}\label{alter}
Theorem $\ref{main}${\rm (i)} or {\rm (iii)} holds if $G_0$ is an alternating group.
\end{lem}

\pf This follows from \cite[Theorem 1.1]{FOS}. \hal

\vspace{4mm}
In view of the previous two lemmas, we can suppose from now on that $G_0$ is a classical simple group. Assume that 
\begin{equation}\label{b8}
b(G) \ge 7.
\end{equation}
We aim to show that conclusion (ii) of Theorem \ref{main} must hold. By the above assumption, the dimension $d\ge 7$, and also  every element of $V^{6}$ is fixed by some element of prime order in  $G\setminus Z$, and so
\begin{equation}\label{v7}
V^{6} = \bigcup_{g \in {\mathcal P}} C_{V^{6}}(g),
\end{equation}
where ${\mathcal P}$ denotes the set of elements of prime order in $G\setminus Z$. 
Now $|C_{V^{6}}(g)| = |C_{V}(g)|^{6}$, and 
\begin{equation}\label{alph}
\dim C_V(g) \le \left\lfloor (1-\frac{1}{\a(g)}) \dim V \right \rfloor,
\end{equation}
 where $\a(g)$ is as defined in the preamble to Lemma \ref{gusa}  (strictly speaking, it is $\a(gZ)$ for $gZ \in G/Z$). Writing $\a = \a(G_0)$, it follows that 
\[
|V|^{6} = q^{6d}  \le |{\mathcal P}|\,q^{6\lfloor d(1-\frac{1}{\a}) \rfloor}.
\]
Since $|G| = |Z|\,|G/Z| \le (q-1)\,|{\rm Aut}(G_0)|$, we therefore have
\begin{equation}\label{crude}
q^{6\lceil d/\a \rceil } \le |{\mathcal P}| < |G| \le (q-1)\,|{\rm Aut}(G_0)|.
\end{equation}

\begin{rem}\label{m24} {\rm Using (\ref{crude}) we can handle the case $G_0=M_{24}$ as follows, completing the proof of Lemma \ref{sporadexcep}: we have $\a(M_{24}) \le 4$ by \cite[2.4]{good}, so (\ref{crude}) yields 
$\frac{6}{4}d < \log_2|M_{24}|$, hence $d\le 18$. By \cite{HM}, this forces $d=11, q=2$, so $G = M_{24} < GL_{11}(2)$. Here $V$ or $V^*$ is a quotient of the binary Golay code of length 24, dimension 12, by a trivial submodule, and we see from \cite[p.94]{atlas} that there is a $G$-orbit on $V$ of size 276 or 759 on which $G$ acts primitively. The base sizes of these actions of $M_{24}$ are less than 7, by \cite{BOW}, and the conclusion follows.}
\end{rem}

Let $q=p^a$, where $p$ is prime. The analysis divides naturally, according to whether or not the underlying characteristic of $G_0$ is equal to $p$ -- that is, whether or not $G_0$ is in the set ${\rm Lie}(p)$.

\begin{lem}\label{liepdash}
Under the above assumption $(\ref{b8})$, $G_0$ is not in ${\rm Lie}(p')$.
\end{lem}

\pf  Suppose $G_0 \in {\rm Lie}(p')$. Lower bounds for $d = \dim V$ are given by \cite{LaS, SZ}, and the values of $\a$ by Lemma \ref{gusa}. Plugging these into (\ref{crude}) (and also using the fact that $d\ge 7$), we see that $G_0$ must be 
one of the following:
\[
\begin{array}{l}
PSp_4(3),\, PSp_4(5),\,Sp_6(2),\,PSp_6(3),\,PSp_8(3), \,  PSp_{10}(3),\\
U_3(3),\, U_4(3),\, U_5(2), \\
\O_7(3),\, \O_8^+(2).
\end{array}
\]
At this point we use \cite{HM}, which gives the dimensions and fields of definition of all the irreducible projective representations of the above groups of dimension up to 250. Combining this information with (\ref{crude}) leaves just the following possibilities:
\[
\begin{array}{|l|l|l|}
\hline
G_0 & d & q \\
\hline
U_5(2) & 10 & 3 \\
U_4(3) & 20 & 2 \\
Sp_6(2) & 7,8& q\le 11 \\
             & 14 & 3 \\
\O_8^+(2) & 8 & q\le 29 \\
\hline
\end{array}
\]
Consider first  $G_0=U_5(2)$. Here $G = \L -I\R \times U_5(2).2 < GL_{10}(3)$, and the Brauer character of this representation of $G$ is given in \cite{atlas}. From this we can read off the dimensions of the fixed point spaces of $3'$-elements of prime order. These are as follows, using Atlas notation:
\[
\begin{array}{r|c|c|c|c|c}
g & 2A,-2A & 2B,-2B & 2C,-2C & 5A & 11AB \\
\hline
\dim C_V(g) & 2,8 & 6,4 & 5,5 & 2& 0 \\
\end{array}
\]
Also $\a\le 5$ by Lemma \ref{gusa}, so (\ref{alph}) gives $\dim C_V(g)\le 8$ for all elements $g\in G$ of order 3. At this point, the inequality 
$|V|^6 \le \sum_{g\in {\mathcal P}} |C_V(g)|^6$ implied by (\ref{v7}) gives 
\[
3^{60} \le |2A|\cdot (3^{12}+3^{48}) +  |2B|\cdot (3^{24}+3^{36}) +  |2C|\cdot (3^{30}+3^{30}) +  |5A|\cdot 3^{12} +  |3ABCDEF|\cdot 3^{48},
\]
where $|2A|$ denotes the size of the conjugacy class of $2A$-elements, and so on. 
This is a contradiction. 

This method works for all the cases in the above table, except  $(G_0,d,q) = (\O_8^+(2),8,3)$; in this case the crude inequality 
$|V|^6 \le \sum_{g\in {\mathcal P}} |C_V(g)|^6$ implied by (\ref{v7}) does not yield a contradiction. Here we have $G \le 2.O_8^+(2) < GL(V) = GL_8(3)$. 
Observe that $O_8^+(2)$ has a subgroup $N = S_3\times O_6^-(2)$, and $N$ is the normalizer of $\langle x\rangle$, where $x$ is an element of order 3.
Then $C_V(x) \ne 0$, and $N$ must fix a 1-space in $C_V(x)$. Moreover, we compute that the minimal base size of 
$O_8^+(2)$ acting on the cosets of $N$ is equal to 4. It follows that there are four 1-spaces in $V$ whose pointwise stabilizer in $G$ is contained in $Z$. Hence 
$b(G) \le 4$ in this case. 
\hal 

\vspace{2mm}
In view of the previous lemmas, from now on we may assume that $G_0 = Cl_n(q_0)$, a classical group over a field $\F_{q_0}$ of characteristic $p$. Recall that $G \le GL(V) = GL_d(q)$ and $G_0 = {\rm soc}(G/Z)$. The next lemma identifies the possible highest weights for $V$ as a module for the quasisimple classical group $E(G)$.

\begin{lem}\label{liep}
Suppose as above that $G_0 = Cl_n(q_0)$, a classical group in ${\rm Lie}(p)$. Then $\F_{q_0}$ is a subfield of $\F_q$, and one of the following holds:
\begin{itemize}
\item[(1)] $V = V(\l)$, where $\l$ is one of the following high weights (listed up to automorphisms of $G_0)$:
\[
\l_1,\,\l_2,\,2\l_1,\,\l_1+p^i\l_1,\,\l_1+p^i\l_{n-1} \,(i>0)
\]
(the last one only for $G_0 = L_n^\e(q_0)$);
\item[(2)] $G_0 = L_n^\e(q_0)\,(n\ge 3)$, $V = V(\l_1+\l_{n-1})$;
\item[(3)] $G_0 = L_n(q_0)\,(7\le n\le 21)$ and $V = V(\l_3)$;
\item[(4)] $G_0 = L_6^\e(q_0)$ and $V = V(\l_3)$;
\item[(5)] $G_0 = L_8^\e(q_0)$ and $V=V(\l_4)$;
\item[(6)] $G_0 = PSp_6(q_0)$ and $V = V(\l_3)$ ($p$ odd);
\item[(7)] $G_0 = PSp_8(q_0)$ and $V = V(\l_3)$ ($p$ odd) or $V(\l_4)$ ($p$ odd);
\item[(8)] $G_0 = PSp_{10}(q_0)$ and $V=V(\l_3)$ ($p=2$);
\item[(9)] $G_0 = P\O_n^\e (q_0)\,(7\le n\le 20,\,n\ne 8)$ and $V$ is a spin module.
\end{itemize}
\end{lem}

\pf Assume first that $q_0>q$. Then by \cite[5.4.6]{KL}, there is an integer $s\ge 2$ such that $q_0 = q^s$ and $d=m^s$, where $m$ is the dimension of an irreducible module for $E(G)$. Note that $m\ge n$ (by the minimal choice of $n$). By (\ref{crude}), 
\[
q^{6m^s/\a} \le (q-1)\,|{\rm Aut}(Cl_n(q^s))|.
\]
 Lemma \ref{gusa} shows that $\a \le n+2$ (excluding the small groups in Lemma \ref{gusa}(vii)), and hence
\[
q^{6m^s/(n+2)} \le (q-1)\,|{\rm Aut}(Cl_n(q^s))| < (q-1)\,q^{s(n^2-1)}\,(2s\log_pq).
\]
Since $m\ge n$, it follows from this that $s=2$ and 
\[
m^2 < \frac{(n+2)(2n^2+1)}{6}.
\]
Now using \cite{Lu}, we deduce that $m=n$ and so
\[
E(G) \le SL_n(q^2) < SL_{n^2}(q).
\]
As in \cite[p.104]{LSbase}, we see that there is a vector $v$ such that $E(G)_v \le SU_n(q)$. By Lemma \ref{sporex}, the base size of an almost simple group with socle $L_n(q^2)$ acting on the cosets of a subgroup containing $U_n(q)$ is at most 4. Hence there are 1-spaces $\d_1,\ldots,\d_4$ whose pointwise stabilizer in $G$ is contained in $Z$, and so $b(G) \le 4$ in this case. This contradicts our initial assumption that $b(G)\ge 7$.

Hence we may assume now that $q_0\le q$, so that $\F_{q_0}$ is a subfield of $\F_q$ by \cite[5.4.6]{KL}. Now (\ref{crude}) gives
\begin{equation}\label{abd}
d < \frac{\a}{6}\left(1+\log_q |{\rm Aut}(G_0)|\right).
\end{equation}
Noting that apart from the case where $G_0 = P\O_8^+(q_0)$, we have $|{\rm Out}(G_0)| \le q$, it now follows using Lemma \ref{gusa} that $d < N$, where $N$ is as defined in Table \ref{Ndef}. In the last row of the table, $\d$ is $\log_q6$ if $G_0 = P\O_8^+(q_0)$, and is 0 otherwise.

\begin{table}[h]
\caption{}\label{Ndef}
\[
\begin{array}{|l|l|}
\hline
G_0 & N \\
\hline
L_n^\e(q_0) & \frac{1}{6}(n+2)(1+n^2),\;n\le 4 \\
&    \frac{1}{6}n(1+n^2),\;n> 4 \\
PSp_n(q_0),\,n\ge 4 & \frac{1}{6}(n+1)\left(2+\frac{1}{2}n(n+1)\right), \;n>4 \\
P\O_n^\e(q_0),\,n\ge 7 & \frac{1}{6}n\left(2+\frac{1}{2}n(n-1)\right)+\d \\
\hline
\end{array}
\]
\end{table}

Now applying the bounds in \cite{Lu} (and also the improved bound for type $A$ in \cite{mart}), we see that with one possible exception, one of the cases (1)-(9) in the conclusion holds. The possible exception is $G_0 = L_4^\e(q_0)$ with $p=3$ and $V=V(\l_1+\l_2)$, of dimension 16. But in this case $G$ does not contain a graph automorphism of $G_0$ (since the weight $\l_1+\l_2$ is not fixed by a graph automorphism), and so \cite[4.1]{GS} implies that we can take $\a=4$ in (\ref{abd}), and this rules out this case. 
\hal

\begin{lem}\label{3to8}
Under the above assumption $(\ref{b8})$, $G_0$ is not as in $(3)-(9)$ of Lemma $\ref{liep}$.
\end{lem}

\pf Suppose $G_0$ is as in (3)--(9) of Lemma \ref{liep}. First we consider the actions of the simple algebraic groups $\bar G$ over $K = \bar \F_q$ corresponding to $G_0$ on the $K\bar G$-modules $\bar V = V \otimes K= V_{\bar G}(\l)$. Define 
\[
M_\l = {\rm min}\left\{{\rm codim}V_\g(g)\, |\, \g \in K^*,\,g \in \bar G\setminus Z(\bar G) \right\}.
\]
By Lemma \ref{lawth}, a lower bound for $M_\l$ is given by ${\rm min}(s_\l,s_{\l'})$, and simple calculations give the following lower bounds:
\[
\begin{array}{|l|l|l|}
\hline
\bar G & \l & M_\l \ge \\
\hline
A_n\,(n\ge 5) & \l_3 & \frac{1}{2}(n-1)(n-2) \\
A_7 & \l_4 & 20 \\
C_3 & \l_3\,(p>2) & 4 \\
C_4 & \l_3\,(p>2) & 12 \\
       & \l_4\,(p>2) & 13 \\
C_5 & \l_3\,(p=2) & 24 \\
D_n\,(n\ge 5) & \l_{n-1},\l_n & 2^{n-3} \\
B_n\,(n\ge 3) & \l_n & 2^{n-2} \\
\hline
\end{array}
\]
Apart from cases (4) and (5) of Lemma \ref{liep}, the group $G/Z$ is contained in $\bar G/Z$; in cases (4) and (5), a graph automorphism of $\bar G$ may also be present.  Thus excluding (4) and (5), we see that (\ref{v7}) gives 
\begin{equation}\label{mbd}
q^{6M_\l} \le |G|.
\end{equation}
The bounds for $M_\l$ in the above table now give a contradiction, except when $\bar G = D_n\,(n\le 6)$ or $B_n\,(n\le 5)$. 

We now consider the cases $\bar G = D_n\,(n\le 6)$ or $B_n\,(n\le 5)$. Since $B_{n-1}(q) < D_n(q) < GL(V)$, it suffices to deal with $\bar G = D_6, D_5$ or $B_3$. 

Suppose $G_0 = D_6^\e(q_0)$ with $\F_{q_0} \subseteq \F_q$. By Lemma \ref{d56}(i), for any element $g\in G$ that is not a scalar multiple of a root element, we have ${\rm codim}C_V(g) \ge 12$; and for root elements $u$, from the above table  we have ${\rm codim}C_V(u) \ge 8$. The number of root elements in $G_0$ is less than $2q^{18}$. Hence (\ref{v7}) gives
\[
|V|^6 = q^{32\times 6} \le 2q^{18}(q-1)\cdot q^{24\times 6} + |G|q^{20\times 6},
\]
which is a contradiction.

Now suppose $G_0 = D_5^\e(q_0)$. We perform a similar calculation, using Lemma \ref{d56}(ii). The number of semisimple elements $s$ of $G$ for which $C_{\bar G}(s)'=A_4$ is at most $|Z|\cdot (q-1)|D_5^\e (q):A_4^\e (q).(q-1)| < 2q^{22}$. The number of root elements in $G_0$ is less than $2q^{14}$, and the number of unipotent elements in the class 
$(A_1^2)^{(1)}$ is less than $2q^{20}$ (these have centralizer in $D_5^\e(q)$ of order $q^{14}|Sp_4(q)|(q-\e)$, see \cite[Table 8.6a]{LSbk}). Moreover, the total number of unipotent elements is at most $q^{40}$. Hence (\ref{v7}) together with Lemma \ref{d56}(ii) gives
\[
q^{16\times 6} \le 2(q^{14}+q^{20})(q-1) q^{12\times 6} + q^{40}(q-1) q^{8\times 6} + 
2q^{22} q^{10\times 6} + |G|q^{8\times 6}.
\]
This is a contradiction.

Next consider $G_0 = B_3(q)$. In the action on the spin module $V$, there is a vector $v$ with stabilizer $G_2(q)$ in 
$B_3(q)$. Hence $b(G)\le 4$ in this case, by Lemma \ref{sporex}(ii).

It remains to handle the cases (4), (5), where $G$ may contain graph automorphisms of $\bar G$. For $G_0 = L_6^\e(q)$ or $L_8^\e(q)$, the conjugacy classes of involutions in the coset of a graph automorphism are given by \cite[\S 19]{AS} for $q$ even, and by \cite[4.5.1]{GLS} for $q$ odd. It follows that the number of such involutions is less than $2q^{21}$ or $2q^{36}$ in case (4) or (5), respectively. For such an involution $g$, by (\ref{alph}) we have $\dim C_V(g) \le 16$ or 60, respectively. All other elements of prime order in $G$ lie in $\bar G\,Z$, hence have fixed point space of codimension at least $M_\l$.  Hence we see that (\ref{v7}) gives
\[
|V|^6 = \left\{\begin{array}{l}
q^{20\times 6} \le |G|\cdot q^{14\times 6} + 2q^{21}\cdot q^{16\times 6}, \hbox{ in case (4)}, \\
q^{70\times 6} \le |G|\cdot q^{50\times 6} + 2q^{36}\cdot q^{60\times 6}, \hbox{ in case (5)}.
\end{array}
\right.
\]
Both of these yield contradictions.

This completes the proof of the lemma. \hal

\begin{lem}\label{case2}
The group $G_0$ is not as in $(2)$ of Lemma $\ref{liep}$.
\end{lem}

\pf  Here $G_0 = L_n^\e(q_0)$ with $n\ge 3$, and $V = V(\l_1+\l_{n-1})$. Suppose first that $\e=+$. Then $G \le PGL_n(q)\,Z$, and  $V$ can be identified with $T/T_0$, where
\[
T = \{A \in M_{n\times n}(q)\,:\,{\rm Tr}(A) = 0\}, \;T_0 = \{\l I_n : n\l = 0\},
\]
and the action of $GL_n(q)$ is by conjugation. By \cite{St}, we can choose $X,Y \in SL_{n-1}(q_0)$ generating $SL_{n-1}(q_0)$. Define 
\[
A = \begin{pmatrix}X& 0 \\ 0 & -{\rm Tr}(X) \end{pmatrix}, \; B = \begin{pmatrix}Y& 0 \\ 0 & -{\rm Tr}(Y) \end{pmatrix}.
\]
Then $GL_n(q)_{A,B} \le \{ {\rm diag}(\l I_{n-1}, \mu)\}$, and hence $b(G) \le 4$.

Now suppose $\e=-$, so that $G \le PGU_n(q)\,Z$, where we take $GU_n(q) = \{g \in GL_n(q^2) : g^Tg^{(q)}=I\}$. Then we can identify $V$ with the $\F_q$-space $S$ modulo scalars, where 
\[
S = \{ A \in M_{n\times n}(q^2)\,:\,{\rm Tr}(A) = 0,\,A^T = A^{(q)} \},
\]
with $GU_n(q)$ acting by conjugation. As in \cite[p.104]{LSbase}, there is a vector $A \in V$ such that $GU_n(q)_A \le N_r$, where $N_r$ is the stabilizer of a non-degenerate $r$-space and $r = \frac{1}{2}n$ or 
$\frac{1}{2}(n-(n,2))$. In the first case, the base size of $PGU_n(q)$ acting on ${\mathcal N}_r$ is at most 5, by Lemma \ref{sporex}(ii) (since in this case $N_r$ is contained in a non-subspace subgroup of type $GU_{n/2}(q) \wr S_2$); and the same holds in the second case, by Theorem \ref{nrbase}. It follows that $b(G)\le 5$, contradicting our assumption (\ref{b8}). \hal

\vspace{2mm}
The proof of Theorem \ref{main} is completed by the following lemma.

\begin{lem}\label{case2}
If $G_0$ is as in $(1)$ of Lemma $\ref{liep}$, then conclusion {\rm (ii)} of Theorem $\ref{main}$ holds.
\end{lem}

\pf Here $G_0 = Cl_n(q_0)$, and $V = V(\l)$ with $\l = \l_1,\,\l_2,\,2\l_1,\,\l_1+p^i\l_1$ or $\l_1+p^i\l_{n-1}$. 

 If $\l = \l_1$, then $d=n$ and $E(G) = Cl_d(q_0)$ is as in part (ii) of Theorem \ref{main}.

Now consider $\l = \l_2$. Here we argue as in the proof of \cite[2.2]{LSbase} (see p.102). If  $V = \wedge^2W$ where $W$ is the natural module for $Cl_n(q_0)$ (with scalars extended to $\F_q$), then the argument provides a vector $v \in V$ such that $SL(W)_v = Sp(W)$, and so application of Lemma \ref{sporex}(ii) gives $b(G) \le b(SL(W)/Sp(W)) \le 5$.
Otherwise, $V$ is equal to $(\wedge^2W)^+$ (which is $f^\perp$ or $f^\perp/\langle f\ra$ in the notation of \cite[p.103]{LSbase}), and the argument gives 
\[
b(G) \le b(Sp_{2k}(q),{\mathcal N}_{r}),
\]
where ${\mathcal N}_{r}$ is the set of non-degenerate subspaces of dimension $r$ and $r = \frac{1}{2}n$ or 
$\frac{1}{2}(n-(n,4))$. As before, Lemma \ref{sporex}(ii) (in the first case) and Theorem \ref{nrbase} (in the second) now  give $b(G)\le 5$.

The case where $\l= 2\l_1$ is similar to the $\l_2$ case, arguing as in \cite[p.103]{LSbase}. Note that $p$ is odd here. If $G_0$ is not an orthogonal group, then $E(G) \le SL(W)$ acting on $V=S^2W$, and there is a vector $v$ such that $SL(W)_v = SO(W)$; hence $b(G) \le b(SL(W)/SO(W)) \le 5$, by Lemma \ref{sporex}(ii). And if $G_0$ is orthogonal, then $V = (S^2W)^+$ (of dimension $\dim S^2W-\d$, $\d\in \{1,2\}$), and we see as in the previous case that $b(G) \le b(O_{2k}(q),{\mathcal N}_{r})$ with 
$r = \frac{1}{2}\left(n-(n,2)\right)$. Hence  Theorem \ref{nrbase} gives $b(G)\le 5$ again. 

Finally, suppose $\l = \l_1+p^i\l_1$ or $\l_1+p^i\l_{n-1}$. Here as in \cite[p.103]{LSbase}, we have $E(G) \le SL(W)=SL_n(q)$ acting on $V = W\otimes W^{(p^i)}$ or $W\otimes (W^*)^{(p^i)}$. We can think of the action of $SL(W)$ on $V$ as the action on $n\times n$ matrices, where $g \in SL(W)$ sends 
\[
A \to g^TAg^{(p^i)} \hbox{ or }g^{-1}Ag^{(p^i)}.
\]
Hence we see that the stabilizer of the identity matrix $I$ is contained in $SU_n(q^{1/2})$ or $SL_n(q^{1/r})$ for some $r>1$, and so as usual Lemma \ref{sporex}(ii) gives $b(G) \le 5$. \hal

\vspace{4mm}
Thsi completes the proof of Theorem \ref{main}.

\end{document}